\documentclass[11pt]{article}
\usepackage{amsmath,amsthm,amsfonts,amssymb}
\usepackage{mathrsfs}
\usepackage{epsfig}
\usepackage{pgf,tikz}
\usepackage{graphicx}
\usepackage{graphics}
\usepackage{setspace}
\usepackage[dvips,arrow,curve]{xy}
\def\bc{\begin{center}}       \def\ec{\end{center}}
\def\ba{\begin{array}}        \def\ea{\end{array}}
\def\be{\begin{equation}}     \def\ee{\end{equation}}
\def\bea{\begin{eqnarray}}    \def\eea{\end{eqnarray}}
\def\beaa{\begin{eqnarray*}}  \def\eeaa{\end{eqnarray*}}

\author{Chunhua Jin, jinch@amss.ac.cn\\Academy of Mathematics and Systems Science, Chinese Academy of Sciences\\Beijing 100190, China\\
Xu-an Zhao, zhaoxa@bnu.edu.cn\\Mathematics Department of Beijing Normal University\\Beijing 100875, China\thanks{The authors are supported by Natural Science Foundation of China (Grant No.11171025)}}

\title{On the Poincar\'{e} series of Kac-Moody Lie algebras}

\begin{document}
\vskip 0.5 true cm

\maketitle
\begin{abstract}
In this paper, we discuss the Poincar\'{e} series of Kac-Moody Lie algebras, especially for indefinite type. Firstly, we compute the Poincar\'{e} series of certain indefinite Kac-Moody Lie algebras whose Cartan matrices have the same type of $2\times 2$ principal sub-matrices. Secondly, we show that the Poincar\'{e} series of Kac-Moody Lie algebras satisfy certain interesting properties. Lastly we give some applications of the Poincar\'{e} series to other fields. Particularly we construct some counter examples to a conjecture of Victor Kac\cite{Kac_85} and a conjecture of Chapavalov, Leites and Stekolshchik\cite{CDR_10}.

\textbf{Keywords:} Kac-Moody Lie algebra, Poincar\'{e} series, Cartan matrix, Flag manifold.
\end{abstract}


\section{Introduction}
To start with, we briefly review some concepts and results about Kac-Moody Lie algebras and their Poincar\'{e} series.\par
 Let $A=(a_{ij})$ be an $n\times n$ integer matrix, $A$ is called a Cartan matrix if it satisfies:

(1) For each $i,a_{ii}=2$;

(2) For $i\not=j,a_{ij}\leq 0$;

(3) If $a_{ij}=0$, then $a_{ji}=0$.

For each Cartan matrix $A$, there is an associated Kac-Moody Lie algebra $g(A)$, modulo its center which is generated by $h_i,e_i,f_i,1\leq i\leq n$ over $\mathbb{C}$,
with the defining relations:

(1) $[h_i,h_j]=0$;

(2) $[e_i,f_i]=h_i,[e_i,f_j]=0,i\not=j$;

(3) $[h_i,e_j]=a_{ij}e_j,[h_i,f_j]=-a_{ij}f_j$;

(4) $\mathrm{ad}(e_i)^{-a_{ij}+1}(e_j)=0$;

(5) $\mathrm{ad}(f_i)^{-a_{i j}+1}(f_j)=0$.

\noindent For details see Kac\cite{Kac_82} and Moody\cite{Moody_68}.

Kac and Peterson\cite{Kac_Peterson_83}\cite{Kac_Peterson_84}\cite{Kac_85} further constructed the Kac-Moody group $G(A)$ with Lie algebra $g(A)$.

Let $M(A)=(m_{ij})_{n\times n}$ be the Coxeter matrix of the Cartan matrix $A$. That is: for $i=j, m_{ij}=1$; for $i\not=j,m_{ij}=2,3,4,6 $ and $\infty$ as $a_{ij}a_{ji}=0,1,2,3 $ and $\geq 4$ respectively, then the Weyl group $W(A)$ of $g(A)$ is: $$W(A)=<\sigma_1,\sigma_2,\cdots,\sigma_n|(\sigma_i\sigma_j)^{m_{ij}}=1, 1\leq i\not=j\leq n >.$$
\indent Each element $w\in W(A)$ has a decomposition $w=\sigma_{i_1}\cdots \sigma_{i_k},1\leq i_1,\cdots,i_k\leq n$. The length of $w$ is defined as the least integer $k$ in all of those decompositions of $w$, denoted by $l(w)$. The Poincar\'{e} series of $g(A)$ is the power series $P(A)=\sum\limits_{w\in W(A)} t^{l(w)}$. Hence $P(A)$ only depends on the structure of Weyl group $W(A)$ and the length function on it.

For the Kac-Moody Lie algebra $g(A)$, there is the Cartan decomposition $g(A)=h\oplus \sum\limits_{\alpha\in \Delta} g_{\alpha}$, where $h$ is the Cartan sub-algebra and $\Delta$ is the
root system. Let $b=h\oplus \sum\limits_{\alpha\in \Delta^+} g_{\alpha}$ be the Borel sub-algebra, then $b$ corresponds to a Borel subgroup $B(A)$ in the Kac-Moody group $G(A)$. The
homogeneous space $F(A)=G(A)/B(A)$ is called the complete flag manifold of $G(A)$. By Kumar\cite{Kumar_02} $F(A)$ is a ind-variety.

The flag manifold $F(A)$ admits a CW-decomposition of Schubert cells which is indexed by the elements of Weyl group $W(A)$. For each $w\in W(A)$, the complex dimension of Schubert variety $X_w$ is $l(w)$. Therefore the Poincar\'{e} series of the flag manifold $F(A)$ defined by Betti numbers is just the Poincar\'{e} series $P(A)$ of $g(A)$. The Poincar\'{e} series $P(A)$ is closely related to the topology of Kac-Moody group $G(A)$ and its flag manifold $F(A)$. For flag manifolds over finite field $F_q,q=p^k$ for a prime number $p$, by the Schubert decomposition, the number of elements in $F(A)$ is just given by the value of Poincar\'{e} series $P(A)$ at $t=q$. This fact can be further used to compute the orders of Kac-Moody group over $F_q$(by using a regularization procedure).

The cohomology of flag manifolds and Poincar\'{e} series of finite type and affine type Kac-Moody Lie algebras are extensively studied by Borel\cite{Borel_53_1}, Bott and Samelson\cite{Bott_Samelson_55}, Bott\cite{Bott_56}, and Chevalley\cite{Chevalley_94}. But for indefinite type, little is known.

For Kac-Moody Lie algebra $g(A)$ of finite type, Poincar\'{e} series has the following form
$$\hspace*{5cm}P(A)=\prod\limits_{i=1}^{n}\dfrac{t^{d_{i}}-1}{t-1} \hspace*{5cm} (1.1)$$
where $d_i$'s are the degrees of basic invariants of $g(A)$.

Let $I$ be a subset of $S=\{1,2\cdots,n\}$, define the principal sub-matrix $A_I=(a_{ij})_{i,j\in I}$. It's obvious that $A_I$ is also a Cartan matrix. Let $W_I(A)$ be the subgroup of $W(A)$ generated by $\{\sigma_i | i\in I\}$, then $W_I(A)$ is the Weyl groups of $g(A_I)$. In the following we also denote $P(A_I)$ by $P_I(A)$.
The Poincar\'{e} series $P(A)$ is also denoted by $P_A(t)$ to emphasize the variable $t$.

According to the classification of Cartan matrix $A=(a_{ij})_{n \times n}$, Kac-Moody algebras are classified into three types:

  (1) $g(A)$ is of finite type if $A$ is positive definite.

  (2) $g(A)$ is of affine type if $A$ is positive semi-definite and has rank $n-1$.

  (3) $g(A)$ is of indefinite type otherwise.

Let $\widetilde A$ be the Extended Cartan matrix of the Cartan matrix $A$ of finite type, then every non-twisted affine Kac-Moody Lie algebra $g(\widetilde A)$ is the central extension of the infinite-dimensional Lie algebra $g(A)\otimes \mathbb{C}[t,t^{-1}]$. Bott\cite{Bott_56} showed that the Poincar\'{e} series is
$$\hspace*{4.3cm}P(\widetilde A)=P(A)\prod\limits_{i=1}^{n}\dfrac{1}{1-t^{d_{i}-1}}\hspace*{4.8cm} (1.2)$$
where $d_i$'s are the degrees of basic invariants of $g(A)$.

For Kac-Moody Lie algebra $g(A)$ of finite type,
$$ \hspace*{5cm}\sum\limits_{I\subset S}(-1)^{|I|}\dfrac{P(A)}{P_{I}(A)}=t^{D(A)} \hspace*{4.8cm}(1.3)$$
where $D(A)$ is the complex dimension of $F(A)$.

For Kac-Moody Lie algebra $g(A)$ of affine or indefinite type,
$$ \hspace*{5cm}\sum\limits_{I\subset S}(-1)^{|I|}\dfrac{P(A)}{P_{I}(A)}=0 \hspace*{5.6cm}(1.4)$$
The Equations (1.3) and (1.4) are given in Steinberg\cite{Steinberg_68} and Humphurays\cite{Humphurays_90}.

\indent Theoretically, one can compute Poincar\'{e} series of any affine and indefinite Kac-Moody Lie algebras through Equation (1.3) and (1.4) by iterations. This reduces the computation of Poincar\'{e} series of Kac-Moody Lie algebras of affine and indefinite types to the finite case.

The paper proceeds as follows: In section 2, we give some concrete computation results of the Poincar\'{e} series of certain Kac-Moody Lie algebras of indefinite type;
In section 3, we derive some properties satisfied by the Poincar\'{e} series. These properties are mainly derived from Equation $(1.4)$ and verify a conjecture given by Gungormez and Karadayi\cite{Gungormez_Karadayi_07}; In section 4, we show the computation of the Poincar\'{e} series for Kac-Moody Lie algebras can be used to compute the Poincar\'{e} series of the generalized flag manifolds of Kac-Moody groups. In section 5, we discuss the application of Poincar\'{e} series to graph theory and define the Poincar\'{e}-Coxeter invariants and the homotopy indices of a graph. In section 6, we construct some counter examples to a conjecture given by Kac\cite{Kac_85}; In section 7, we derive the condition needed for the conjecture of Chapavalov, Leites and Stekolshchik\cite{CDR_10} to be true and give some counter examples.
\section{Computation of the Poincar\'{e} series for certain Kac-Moody Lie algebras}

\hspace{.9cm}From the definition of Weyl group and Poincar\'{e} series, we see that Poincar\'{e} series are determined only by the products $a_{ij}a_{ji}$, for all $i$, $j$. If we change one $a_{ij}a_{ji}>4$ to $4$, then the Coxeter matrix is unchanged. Therefore we have

\noindent\textbf{Proposition 1.}
{\sl Let $A=(a_{ij})_{n\times n}$ be a Cartan matrix, define a Cartan matrix $A'=(a'_{ij})_{n\times n}$. For $i>j$, $(a'_{ij},a'_{ji})$ is given by the product $a_{ij}a_{ji}$ as follows,
$$\begin{tabular}{|c|c|c|c|c|c|}
  \hline
  $a_{ij}a_{ji}$ & $0$ & $1$ & $2$ & $3$ & $\geq 4$ \\\hline
  $(a'_{ij},a'_{ji})$ & $(0,0)$ & $(-1,-1)$ & $(-2,-1)$ & $(-3,-1)$ & $(-4,-1)$ \\
  \hline
\end{tabular}$$\\
then $P(A)=P(A')$}.

Therefore for the computation of Poincar\'{e} series, we only need to consider Cartan matrix of  form:
$${\footnotesize A=\left(\begin{array}{cccc}
        2 & l(a_{21}) & \cdots & l(a_{n1}) \\
         a_{21} & 2 & \cdots & l(a_{n2}) \\
          \vdots &   \vdots & \ddots &\vdots\\
         a_{n1} & a_{n2} & \cdots & 2 \\
       \end{array}
     \right)}~~(\ast),$$
where $-4\leq a_{ij}\leq 0$ and $l(x)= \left\{\begin{array}{ll}
                                       0,& x=0;\\
                                        -1, & \hbox{otherwise.}
                                        \end{array}
                                        \right. $

\noindent\textbf{Theorem 1.} {\sl Let $P_n$ be the Poincar\'{e} series of an $n\times n$ Cartan matrix $A$ in which all order 2 principal sub-matrices having the same type, and $P_2$ denotes the Poincar\'{e} series of the order 2 principal sub-matrix, then
$$P_n=\dfrac{(1+t)P_2}{\frac{(n-1)(n-2)}{2}tP_2-\frac{(n+1)(n-2)}{2}P_2+\frac{n(n-1)}{2}(1+t)}$$}
\noindent \textbf{Proof:} we prove the theorem by induction on $n$.

(i) For $n=2$, Theorem 1 can be checked directly.

(ii) Suppose Theorem 1 is true for $k<n$, by Equation $(1.4)$, we get
$$1-{n \choose 1}\frac{1}{P_1}+{n \choose 2}\frac{1}{P_2}-{n \choose 3}\frac{1}{P_3}+\cdots +(-1)^{n-1}{n \choose n-1}\frac{1}{P_{n-1}}=\frac{(-1)^{n-1}}{P_n}$$
By $P_1=1+t.$ and the induction assumption, we have $$\sum\limits_{i=0}^{n-1}(-1)^{i}{n \choose i}\frac{\frac{(i-1)(i-2)}{2}tP_2-\frac{(i+1)(i-2)}{2}P_2+\frac{i(i-1)}{2}(1+t)}{(1+t)P_2}=\frac{(-1)^{n-1}}{P_n}$$
Since $$\sum\limits_{i=0}^{n-1}(-1)^i{n \choose i}i=(-1)^{n-1}n\text{ and }\sum\limits_{i=0}^{n-1}(-1)^i{n \choose i}i^2=(-1)^{n-1}n^2.$$
We can show $$(-1)^{n-1}\frac{(1+t)P_2}{\frac{(n-1)(n-2)}{2}tP_2-\frac{(n+1)(n-2)}{2}P_2+\frac{n(n-1)}{2}(1+t)}=\frac{(-1)^{n-1}}{P_n}.$$
So the theorem holds for $k=n$.\qed

By Theorem 1, we get\\
\noindent\textbf{Corollary 1.} {\sl Let $A=(a_{ij})_{n\times n}$ be a Cartan matrix as in Theorem 1, then\\

\noindent $P(A)=\left\{
             \begin{array}{cl}
               \dfrac{P(A_2)}{\frac{(n-1)(n-2)}{2}t^{3}-(n-2)t^{2}-(n-2)t+1}& \hbox{if}\quad \forall i\neq j, a_{ij}a_{ji}=1;\\
               \dfrac{P(B_2)}{\frac{(n-1)(n-2)}{2}t^{4}-(n-2)t^{3}-(n-2)t^{2}-(n-2)t+1} &\hbox{if}\quad\forall i\neq j, a_{ij}a_{ji}=2;\\
               \dfrac{P(G_2)}{\frac{(n-1)(n-2)}{2}t^{6}-(n-2)(t^5+t^4+t^3+t^2+t)+1} & \hbox{if}\quad\forall i\neq j, a_{ij}a_{ji}=3;\\
               \dfrac{1+t}{1-(n-1)t} & \hbox{if}\quad\forall i\neq j, a_{ij}a_{ji}\geq 4.
         \end{array}
           \right.$\vspace*{2mm}}\par
The computation in the proof of Theorem 1 is just a sample.
The same method can be used to give a complete computation results for all the Poincar\'{e} series of finite and affine types(untwisted case and twisted case). For various special kinds of Cartan matrices of indefinite type, the method is also useful to get the various kinds of formulas for Poincar\'{e} series. But there is not a single simple formula, so we don't discuss it further.

\section{Some properties of the Poincar\'{e} series}

It is apparent that two different Cartan matrices may have the same Poincar\'{e} series. The following theorem shows that for Cartan  matrix of form $(\ast)$ without element 0, the corresponding Poincar\'{e} series is determined by the set $S(A):=\{a_{ij}| 1\leq j <i \leq n\}$(with multiplicity).

\noindent\textbf{Theorem 2.} {\sl For a Cartan matrix  $A$ of form $(\ast)$ with $ a_{ij}\neq 0,\hspace*{0.1cm} \forall i, j$, the corresponding Poincar\'{e} series is invariant under exchanging of any two elements $a_{ij}$ and $a_{i'j'}$ below diagonal. }

\noindent\textbf{Proof:} set $$g(a):=\left\{
             \begin{array}{ll}
               \displaystyle{\frac{(1-t^{2})(1-t^{3})}{(1-t)^{2}}}, & a=-1;\vspace*{1.414mm} \\
               \displaystyle{\frac{(1-t^{2})(1-t^{4})}{(1-t)^{2}}}, & a=-2;\vspace*{1.414mm}\\
               \displaystyle{\frac{(1-t^{2})(1-t^{6})}{(1-t)^{2}}}, & a=-3;\vspace*{1.414mm} \\
               \displaystyle{\frac{1+t}{1-t}}, & a=-4.
             \end{array}
           \right.$$
In fact, $g(a)$ is the Poincar\'{e} series of Cartan matrix ${\footnotesize \left(
                                                                    \begin{array}{rr}
                                                                      2 & -1 \\
                                                                      a & 2 \\
                                                                    \end{array}
                                                                  \right)}$.\par

To prove the theorem, we only need to prove

\noindent {\bf Assertion:} For a Cartan matrix $A=(a_{ij})_{n\times n}$ without element 0, $P(A)$ is a symmetric rational function of $\{g(a_{ij})|1\leq j<i \leq n\}$.\par
  Now, we start to prove this assertion. \par

Now we prove the assertion by induction on $n$. Notice that, if there is no element 0 in $A$, then all principal sub-matrices with order $\geq 3$ in $A$ are of indefinite type.\par\par
(i)For $n=3$, suppose $A=(a_{ij})_{3\times 3}$ without element 0, by Equation $(1.4)$, we have
$$P(A)-3\dfrac{P(A)}{1+t}+\dfrac{P(A)}{g(a_{21})}+\dfrac{P(A)}{g(a_{31})}+\dfrac{P(A)}{g(a_{32})}=1.$$
hence
$$P(A)=\dfrac{(1+t)g(a_{21})g(a_{31})g(a_{32})}{(t-2)g(a_{21})g(a_{31})g(a_{32})+(t+1)(g(a_{21})g(a_{31})
+g(a_{21})g(a_{32})+g(a_{31})g(a_{32}))}.$$
So $P(A)$ is a symmetric rational function of $g(a_{21}), g(a_{31}), g(a_{32})$.

(ii)Suppose the assertion hold for all $k$ with $3\leq k < n$, and $A=(a_{ij})_{n\times n}$ be a Cartan matrix without element 0.
By the induction assumption, for any $I\subsetneqq S$, $P_I(A)$ is a symmetric rational function of $\{g(a_{ij})|i,j\in I ,j < i\}$. We define an action of $S_{\frac{n(n-1)}{2}}$ (permutation group of $\frac{n(n-1)}{2}$-elements) on the set of rational functions $g(a_{ij})$'s by the action of $S_n$ on index $(i,j)$. Then we get an action of $S_{\frac{n(n-1)}{2}}$ on the rational functions of $g(a_{ij})$'s.

Note the set $$\{\ P_{I}(A)\ |\ |I|=k, I\subsetneqq S\},k<n$$ is invariant under the $S_{\frac{n(n-1)}{2}}$ action.
Therefore, for any $k<n$, $\sum_{|I|=k, I\subsetneqq S}\dfrac{1}{P_{I}(A)}$ is a symmetric rational function of $\{g(a_{ij})|1\leq j<i \leq n\}$. And by equation $$(-1)^{|S|}+\sum\limits_{I\subsetneqq  S}(-1)^{|I|}\frac{P(A)}{P_{I}(A)}=\sum\limits_{I\subset S}(-1)^{|I|}\frac{P(A)}{P_{I}(A)}=0$$
we get that $P(A)$ is a symmetric rational function of $\{g(a_{ij})|1\leq j<i \leq n\}$. This finishes the proof. \qed

The following example shows that the condition \lq\lq there is no element 0 in $A$\rq\rq   can't be discarded.

\noindent\textbf{Example 1:} For Cartan matrices$${ A=\left(
          \begin{array}{rrrr}
            2 & -1 & -1 & -1 \\
            -1 & 2 & 0 & -1 \\
            -2 & 0 & 2 & -1 \\
            -3 & -3 & -1 & 2 \\
          \end{array}
        \right)}\quad \text{and} \quad {B=\left(
          \begin{array}{rrrr}
            2 & 0 & -1 & -1 \\
            0 & 2 & -1 & -1 \\
            -1 & -1 & 2 & -1 \\
            -2 & -3 & -3 & 2 \\
          \end{array}\right)}$$
\noindent their Poincar\'{e} series are
$$P(A)={\footnotesize \dfrac{(t^3+1)(t^3+t^2+t+1)(t^2+t+1)(t+1)}{2t^9+t^8-3t^5-3t^4-2t^3-2t^2-t+1}}$$
and
$$P(B)={\footnotesize\dfrac{(t^5+t^4+t^3+t^2+t+1)(t^3+t^2+t+1)}{2t^8-t^7+t^6-2t^5-t^4-2t^3-2t+1}}.$$ Therefore even if the sets $S(A)=\{a_{ij}|1\leq j<i\leq n\}$ and $S(B)=\{b_{ij}|1\leq j<i\leq n\}$ (with multiplicity) are same, the Poincar\'{e} series are different because $A$ and $B$ contain element 0. If we replace $(0, 0)$ in $A$, $B$ by $(-4, -1)$, we get the new Cartan matrices
$$A'={ \left(
          \begin{array}{rrrr}
            2 & -1 & -1 & -1 \\
            -1 & 2 & -1 & -1 \\
            -2 & -4 & 2 & -1 \\
            -3 & -3 & -1 & 2 \\
          \end{array}
        \right)}\text{  and  } B'={ \left(
          \begin{array}{rrrr}
            2 & -1 & -1 & -1 \\
            -4 & 2 & -1 & -1 \\
            -1 & -1 & 2 & -1 \\
            -2 & -3 & -3 & 2 \\
          \end{array}\right)},$$
          and their Poincar\'{e} series are
$$P(A')=P(B')=\dfrac{(t^5+t^4+t^3+t^2+t+1)(t^3+t^2+t+1)}{2t^8-2t^7-2t^5-3t^4-2t^3-t^2-2t+1}$$

Due to the symmetry in Theorem 2, we see that the number of Poincar\'{e} series corresponding to all Cartan matrices without element 0 is much smaller than the number it seems to be. Denote the number of Poincar\'{e} series corresponding to order $n$ Cartan matrices without element 0 by $K(n)$, then through Theorem 2 and some Combinatorial computation, we get $K(n)\leq {\frac{n(n-1)}{2}+3 \choose 3}$. This motivates us to make the following

\noindent {\bf Conjecture 1:} $K(n)={\frac{n(n-1)}{2}+3 \choose 3}$.

The conjecture is verified in cases $n=3,4,5,6$ by figuring out all the $K(n)$ Poincar\'{e} series.

Besides the above, there is another interesting property. It's shown in Gungormez and Karadayi\cite{Gungormez_Karadayi_07} that for Poincar\'{e} series of hyperbolic Kac-Moody Lie algebras, their explicit forms seem to be the ratio of the Poincar\'{e} series of a properly
chosen finite Lie algebra and a polynomial of finite degree. In the following, we will see this result can be extended to all Kac-Moody algebras.

\noindent\textbf{Theorem 3.}
{\sl For any Kac-Moody Lie algebras  $g$ with Cartan matrix $A$, the Poincar\'{e} series can be written as a rational function whose numerator is the least common multiple of those Poincar\'{e} series of all finite type principal sub-matrices of $A$.}

\noindent {\bf Proof: } by Equation $(1.4)$, we have $$\frac{(-1)^{n+1}}{P(A)}=\sum\limits_{I\subsetneqq S } \frac{(-1)^{|I|}}{P_{I}(A)}. $$
\indent For those $P_I(A)$'s with affine or indefinite Cartan matrices $A_I$, substitute the item $\displaystyle{\frac{(-1)^I}{P_I(A)}}$ by the same formula for $A_I$, do the same operation successively until the denominator of each item of the right side of the resulted formula is the Poincar\'{e} series of a finite type sub-matrix.  This proves the theorem. \qed

This shows that the Poincar\'{e} series of a Kac-Moody Lie algebra is a rational function, and its numerator can be chosen to be the least common multiple of the Poincar\'{e} series of the finite type principal sub-matrix of $A$. By multiplying suitable polynomial factor if needed, the numerator can be the Poincar\'{e} series of a properly chosen finite type Lie algebra. Thus we have shown that for an infinite dimensional Kac-Moody Lie algebra $g(A)$, its Poincar\'{e} series can be written as the ratio of the Poincar\'{e} series $P(A')$ of a properly
chosen finite type Lie algebra ${g(A')}$ and a denominator polynomial $Q'$ of finite degree.
                       $$P(A)=\frac{P(A')}{Q'}$$
This proves a conjecture given by Gungormez and Karadayi\cite{Gungormez_Karadayi_07}.

\noindent\textbf{Example 2:} Take Cartan matrix  { $A=\left(
                                   \begin{array}{rrrr}
                                     2 & -1 & -1 & 0 \\
                                     -1 & 2 & -1 & 0 \\
                                     -1 & -1 & 2 & -1 \\
                                     0 & 0 & -3 & 2 \\
                                   \end{array}
                                 \right)$}, its Coxeter graph is
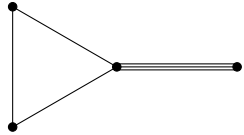
\begin{center}
\begin{tikzpicture}[scale=0.8]
{\filldraw [black] (10,0.5) circle (2pt);}
{\filldraw [black] (12,0.5) circle (2pt);}
{\filldraw [black] (10-1.73,-0.5) circle (2pt);}
{\filldraw [black] (10-1.73,1.5) circle (2pt);}

{\draw (10,0.5) -- +(-1.73, 1);}
{\draw (10,0.5) -- +(-1.73, -1);}

{\draw (10,0.5) -- +(2,0);}
{\draw (10,0.55) -- +(2,0);}
{\draw (10,0.45) -- +(2,0);}

{\draw (10-1.73,-0.5) --(10-1.73,1.5);}

\end{tikzpicture}
\end{center}

\noindent For $A$ we see that all its finite type principal sub-matrices are $A_1$, $A_2$, $G_2$, $A_2\oplus A_1$, and $A_1\oplus A_1$. Since the Poincar\'{e} series of $G_2$ is the least common multiple the of Poincar\'{e} series of these Cartan matrices, the numerator of Poincar\'{e} series of $A$ is just the Poincar\'{e} series of $G_2$.

By calculation, Poincar\'{e} series of $A$ is:
$$\frac{(1+t)(1+t+t^2+t^3+t^4+t^5)}{2t^6-t^5-t^4+t^3-2t+1}$$

\noindent\textbf{Example 3:} For Cartan matrix corresponding to the Dynkin diagram as follows.

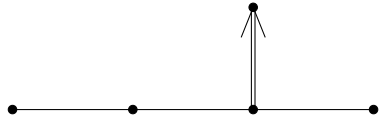
\begin{center}
\begin{tikzpicture}[scale=0.8]
\foreach \x in {10,12,14,16}
\foreach \y in {0.5}
{\filldraw [black] (\x,\y) circle (2pt);}
{\filldraw [black] (14,2.2) circle (2pt);}

\foreach \x in {10,12,14}
\foreach \y in {0.5}
{\draw (\x,\y) -- +(2,0);}

{\draw (14,2.2) -- +(0.2,-0.5);}
{\draw (14,2.2) -- +(-0.2,-0.5);}

{\draw (13.97,0.5) -- +(0,1.7);}
{\draw (14.03,0.5) -- +(0,1.7);}
\end{tikzpicture}
\end{center}

We see that all its finite principal sub-matrices are of type $A_1,\,A_2,\,A_3,\,A_4,\,B_2,\,B_3,\,B_4$, $A_1\oplus B_3,\,A_2\oplus A_1\oplus A_1,\,A_1\oplus A_1\oplus A_1,\,A_1\oplus B_2,\,A_2\oplus A_1$ and $A_1 \oplus A_1$. The least common multiple of their Poincar\'{e} series  is \footnotesize{$\prod\limits_{i=2, 4, 6, 8, 5}\dfrac{t^{i}-1}{t-1}$}, \normalsize{}it is the Poincar\'{e} series of $D_5$. By calculation, Poincar\'{e} series of above Dynkin diagram is
{\footnotesize $$ \dfrac{(t^7+t^6+t^5+t^4+t^3+t^2+t+1)(t^5+t^4+t^3+t^2+t+1)(t^4+t^3+t^2+t+1)(t^3+t^2+t+1)(t+1)}
{-t^{19}-t^{18}-t^{17}-t^{16}+t^{14}+t^{13}+2t^{12}+2t^{11}+2t^{10}+t^9+t^8-t^7-t^6-2t^5-t^4-t^3+1}$$}

By the way, for case 1, 2 and 3 in Corollary 1, the numerators of Poincar\'{e} series are $P(A_{2}),P(B_{2})$ and $P(G_{2})$ respectively. For Cartan matrices in case 1, 2 and 3, their finite type principal sub-matrices are of type $A_{2},\,B_{2}$ and $G_{2}$ respectively. And for case 4, $1+t$ is just the Poincar\'{e} series of $A_1$, the only finite type principal sub-matrix.

\section{The Poincar\'{e} series of generalized flag manifolds}

For each $I\subset S$, there is a parabolic subgroup $G_I(A)$ corresponding to the Cartan matrix $A_I$. The inclusions $B(A)\subset G_I(A)\subset G(A)$ induce a fibration $$G_I(A)/B(A)\stackrel{i}{\longrightarrow} G(A)/B(A)\stackrel{\pi}{\longrightarrow} G(A)/G_I(A). $$
$G_I(A)/B(A)\cong F(A_I)$, and $G(A)/G_I(A)$ is called a generalized flag manifold. This gives a fibration $F(A_I)\stackrel{i}{\longrightarrow} F(A)\stackrel{\pi}{\longrightarrow} G(A)/G_I(A)$. The generalized flag manifold $G(A)/G_I(A)$ also has a decomposition of Schubert cells which are indexed by the coset $W(A)/W_I(A)$. By Bernstein, Gel'fand and Gel'fand\cite{BGG_73}, $W(A)/W_I(A)$ can be regarded as subset $$W^I(A)=\{w\in W|\text{for all } w'\in W_I(A), l(ww')\geq l(w)\}.$$ of $W(A)$. The homomorphism $i^*:H^*(F(A))\to H^*(F(A_I))$ and $\pi^*:H^*(G(A)/G_I(A))\to H^*(F(A))$ have good properties.
See \cite{BGG_73} for details.

\noindent{\bf Lemma 1: }1. $W(A)=W^I(A)\cdot W_I(A)$ and for each $w\in W^I(A),w'\in W_I(A), l(ww')=l(w)+l(w')$.

2. Schubert variety $X_w, w\in W_I(A)$ forms an additive basis of $H^*(F(A_I))$ and $i^*$ is surjective. In fact, for each $w\in W(A)$, if $w\in W_I(A)$, then $i^*(X_w)=X_w$; otherwise $i^*(X_w)=0$. Here for $w\in W_I(A)$, we identify $w\in W_I(A)$ with its image in $W(A)$.

3. Schubert variety $X_w, w\in W^I(A)$ forms an additive basis of $H^*(G(A)/G_I(A))$ and $\pi^*$ is injective. In fact for each $w\in W^I(A)$, $\pi^*(X_w)=X_w$; Here we identify $w$ with $wW_I(A)$ for $w\in W^I(A)$.

As a consequence of the Lemma 1, we get

\noindent{\bf Corollary 2:} The Leray-Serre spectral sequence of fibration $F(A_I)\stackrel{i}{\longrightarrow} F(A)\stackrel{\pi}{\longrightarrow} G(A)/G_I(A)$ collapses at $E_2$-item and $P(A)=P_I(A) \cdot P(G(A)/G_I(A))$.

Hence if we have computed the Poincar\'{e} series of $P(A)$ and $P_I(A)$, the Poincar\'{e} series of the generalized flag manifolds $G(A)/G_I(A)$ is $\displaystyle{\frac{P(A)}{P_I(A)}}$.

\section{Application to graph theory}
The Poincar\'{e} series of Kac-Moody Lie algebras can be used to construct invariants of graphs.

\noindent{\bf Definition 1:} Let $\Gamma$ be an undirected graph without selfloop, suppose the set of vertices of $\Gamma$ to be $\{1,2,\cdots,n\}$ and the vertices are connected by edges with multiplicity. Let $M(\Gamma)$ be the Coxeter matrix $(m_{ij})_{n\times n}$ defined by $m_{ij}=1$ if $i=j$; $m_{ij}=2,3,4,6 $ and $\infty$ if the multiplicity of edge between vertices $i,j$ is $0,1,2,3 $ and $\geq 4$ respectively. The Coxeter group $W(\Gamma)$ of graph $\Gamma$ is $$W(\Gamma)=<\sigma_1,\cdots,\sigma_n|(\sigma_i\sigma_j)^{m_{ij}}=1, 1\leq i,j\leq n>.$$
The Poincar\'{e}-Coxeter invariant of $\Gamma$ is defined as the Poincar\'{e} series of $W(A)$, denoted by $P(\Gamma)$.

According to the definition, the properties in the previous sections can be interpreted naturally as properties of the Poincar\'{e}-Coxeter invariants.
For example, $P(\Gamma)$ is a rational function of variable $t$; If the multiplicity of an edge in $\Gamma$ is great than $4$, then changing the multiplicity to $4$ doesn't alter the invariant; If any two vertices of $\Gamma$ are connected by edges, then exchange of the multiplicities of edges doesn't alter the invariant;






The following definition is motivated by another work of the authors.

\noindent{\bf Definition 2:} Let $\Gamma$ be a graph as above, then the homotopy indices $i_1,i_2,\cdots,i_k,\cdots$ are integers satisfying $\displaystyle{P(\Gamma)=\frac{1}{(1-t)^n}\prod\limits_{k=1}^\infty (1-t^{2k-1})^{i_{2k-1}}(1-t^{2k})^{-i_{2k}}}$.

By results of the authors\cite{Zhao_Jin_Zhang_12}, the sequences $i_k,k>0$ is well defined and we have

The homotopy indices of the Coxeter graph of Kac-Moody groups of finite and affine types are tabulated as follows.

\begin{tabular}{|l|l|l|l|}
  \hline
  $\Gamma$ & nonzero homotopy indices & $\Gamma$ & nonzero homotopy indices \\  \hline
  $A_n$ & $i_k=1,2\leq k\leq n+1 $& $\widetilde A_n$ &$i_k=1,k=n+1  $\\  \hline
  $B_n$ & $i_k=1, k=2,4,\cdots,2n$& $\widetilde B_n$ &$ i_k=1, k=2,3,4,\cdots,2n-1,2n; $\\ \hline
  $D_{2n}$ & $i_k=1, k=2,4,\cdots,$& $\widetilde D_{2n}$ &$i_k=1, k=2,3,4,\cdots,2n-2; $\\
  $n\geq 2$ & $4n-4,4n-2,2n $& $n\geq 2 $ &$ i_k=1, k=2n,2n+1,\cdots,4n-2;$\\
  $$ & $$& $ $ &$ i_{2n-1}=2 $\\ \hline
  $D_{2n+1}$ & $i_k=1, k=2,4,\cdots,$& $\widetilde D_{2n+1} $ &$ i_k=1, k=2,3,4,\cdots,2n-1;$\\
  $n\geq 2$ & $4n-2,4n,2n+1$& $n\geq 2 $ &$ i_k=1, k=2n+1,2n+2,\cdots,4n; $\\ \hline
  $G_2$ & $i_k=1, k=2,6$& $\widetilde E_2$ &$ i_k=1, k=2,5,6$\\  \hline
  $F_4$ & $i_k=1, k=2,6,8,12 $& $\widetilde F_4$ &$i_k=1, k=2,5,6,7,8,11,12;  $\\ \hline
  $E_6$ & $i_k=1, k=2,5,6,8,9,12  $& $\widetilde E_6$ &$i_k=1, k=2,4,6,7,9,11,12;   $\\ \hline
  $E_7$ & $i_k=1, k=2,6,8,10, $& $\widetilde E_7$ &$i_k=1, k=2,5,6,7,8,9,10,11,12,13,  $\\
  $$ & $12,14,18 $& $ $ &$ 14,17, 18$\\ \hline
  $E_8$ & $i_k=1, k=2,8,12,14,18, $& $\widetilde E_8$ & $ i_k=1, k=2,7,8,11,12,13,14,17 ,18, $\\
  $$ & $20,24,30$& $ $ &$ 19,20,23,24,29,30$\\ \hline

 \end{tabular}

\noindent{\bf Proposition 2:} The homotopy indices of $\Gamma$ contain the same amount of information as the Poincar\'{e}-Coxeter invariant of graph $\Gamma$.

The homotopy indices of a Coxeter graph is important both for the research of combinatoric graph theory and for the research of algebraic topology of Kac-Moody groups.
For example $i_2$ is determined by the number of loops in $\Gamma$, see \cite{Kac_85}. For more connection between the homotopy indices and the rational homotopy groups of Kac-Moody groups, see \cite{Zhao_Jin_Zhang_12}.

So a natural question is

\noindent{\bf Question 1}: Determine the homotopy indices $i_1,i_2,\cdots, i_k,\cdots$ from the information of the graph $\Gamma$.










\noindent{\bf Example 4:} For the graph $\Gamma$ with $3$ vertices and edges with multiplicities $2,3,4$, the Poincar\'{e}-Coxeter invariant is
$$P(\Gamma)=\frac{(1-t^4)(1-t^6)}{(1-t)^2(-t^7-t^6 -2t^5- t^4 - 2t^3 -t+1)}.$$
\indent Let $g$ be the graded free Lie algebra over $\mathbb{C}$ with generators $x_1,x_3,y_3,x_4,x_5,y_5,x_6,x_7$(set $\deg x_i=\deg y_i=i$) and $d_k$ be the dimension of the degree $k$ homogeneous component of $g$. The homotopy indices of graph $\Gamma$ be $i_1=0,i_4=d_4-1,i_6=d_6-1$ and the other $i_k=d_k$.

\section{Disproof of a conjecture of Victor Kac}
In Kac\cite{Kac_85}, page 204, Victor Kac gave a conjecture on the Poincar\'{e} series of Kac-Moody Lie algebras of indefinite type.

To state the conjecture, we need the following definition.

\noindent{\bf Definition 3:} An $n\times n$ Cartan matrix $A$ is called symmetrizable if there exists an invertible diagonal matrix $D$ and a symmetric matrix $B$
such that $A=DB$. $g(A)$ is called a symmetrizable Kac-Moody Lie algebra if its Cartan matrix is symmetrizable.

Put $\epsilon=1$ or $0$ according to $A$ is symmetrizable or not and $C(t)=P(A)(1-t)^n(1-t^2)^{-\epsilon}$, Kac gave the following conjecture.

\noindent{\bf Conjecture 2: } $\displaystyle{C(t)=\frac{1}{1-B(t)}}$, where $B(t)=b_2t^2+b_3 t^2+\cdots+$ and $b_i\geq 0$.

We construct some counter examples to Kac's conjecture.

\noindent{\bf Example 5: }For the Cartan matrices
$${A=
\left(
    \begin{array}{rrrrrr}
      2 & -1 & 0 & 0 & -1 & 0 \\
      -1 & 2 & -1 & 0 & 0 & 0 \\
      0 & -1 & 2 & -1 & 0 & 0 \\
      0 & 0 & -1 & 2 & -1 & 0 \\
     -1 & 0 & 0 & -1 & 2 & -1 \\
      0 & 0 & 0 & 0 & -1 & 2 \\
    \end{array}
  \right)
} \text{ and } {B=
\left(
    \begin{array}{rrrrrr}
      2 & -1 & 0 & 0 & 0 & 0 \\
      -1 & 2 & -1 & 0 & -1 & -1 \\
      0 & -1 & 2 & -1 & 0 & 0 \\
      0 & 0 & -1 & 2 & 0 & 0 \\
      0 & -1 & 0 & 0 & 2 & 0 \\
      0 & -1 & 0 & 0 & 0 & 2 \\
    \end{array}
  \right)
}$$
with Coxeter graphs as follows.
\begin{center}
\begin{tikzpicture}[scale=0.8]
\foreach \x in {10,12,14,16}
\foreach \y in {0.5}
{\filldraw [black] (\x,\y) circle (2pt);}
{\filldraw [black] (11,-1.23) circle (2pt);}
{\filldraw [black] (13,-1.23) circle (2pt);}

\foreach \x in {10,12,14}
\foreach \y in {0.5}
{\draw (\x,\y) -- +(2,0);}
{\draw (10,0.5) -- +(1,-1.73);}
{\draw (11,-1.23) -- +(2,0);}
{\draw (13,-1.23) -- +(1,1.73);}

\foreach \x in {10,12,14,16}
\foreach \y in {0.5}
{\filldraw [black] (\x+10,\y) circle (2pt);}

{\filldraw [black] (12+10,-1.5) circle (2pt);}
{\filldraw [black] (12+10,2.5) circle (2pt);}

\foreach \x in {10,12,14}
\foreach \y in {0.5}
{\draw (\x+10,\y) -- +(2,0);}
{\draw (12+10,0.5) -- +(0,2);}
{\draw (12+10,0.5) -- +(0,-2);}

\draw (22.5,-2)  node{Coxeter graph of $A$};
\draw (12.5,-2)  node{Coxeter graph of $B$ };
\end{tikzpicture}
\end{center}
The Poincar\'{e} series are
$$P(A)=\frac{(1-t^2)(1-t^5)(1-t^6)(1-t^8)}{(1-t)^4(t^{16}-t^{15}-t^{11}+t^9-t^8+t^7+t^3-2t+1)}.$$
and
$$P(B)=\frac{(1-t^2)^2(1-t^4)(1-t^6)(1-t^8)}{(1-t)^5(t^{16}-t^{14}-t^{12}+t^8-t^7+2t^5+t^4-t^3-t^2-t+1)}.$$

$A$ is a symmetric matrix, so
$$C(t)=P(A)\cdot(1-t)^6\cdot(1-t^2)^{-1}$$
From $\displaystyle{C(t)=\frac{1}{1-B(t)}}$, $B(t)=1-C(t)^{-1}$, the Taylor expansion of
 $$B(t)=1-\frac{(1-t)^4(t^{16}-t^{15}-t^{11}+t^9-t^8+t^7+t^3-2t+1)}{(1-t^2)(1-t^5)(1-t^6)(1-t^8)}\cdot \frac{1-t^2}{(1-t)^6}.$$
at $t=0$ is:
$$t^2+t^3+t^4+t^7+t^8+t^9-t^{15}-2t^{16}-2t^{17}-3t^{18}-3t^{19}+O(t^{20})$$
For Cartan matrix $B$, the Taylor expansion of
$$B(t)=1-\frac{(1-t)^5(t^{16}-t^{14}-t^{12}+t^8-t^7+2t^5+t^4-t^3-t^2-t+1)}{(1-t^2)^2(1-t^4)(1-t^6)(1-t^8)}
\cdot \frac{1-t^2}{(1-t)^6}$$
at $t=0$ is:
$$2t^3+t^5-t^6+3t^7-2t^8+3t^9-2t^{10}+5t^{11}-4t^{12}+4t^{13}-5t^{14}+8t^{15}-8t^{16}+6t^{17}-10t^{18}+10t^{19}+O(t^{20})$$
\section{On the conjecture of Chapavalov, Leites and Stekolshchik}
In their paper\cite{CDR_10}, page 208, Chapavalov, Leites and Stekolshchik made a conjecture on the Poincar\'{e} series of infinite dimensional Kac-Moody Lie algebras $g(A)$.

For a rational function $P(t)$, write $P(t)$ as form $P(t)=\displaystyle{\frac{Q(t)}{R(t)}},$ define $\deg P(t)= \deg Q(t)- \deg R(t)$. It doesn't depend on how $Q(t)$ and $R(t)$ are chosen.

\noindent{\bf Conjecture 3: }(Chapavalov-Leites-Stekolshchik) The Poincar\'{e} series of an infinite dimensional Kac-Moody Lie algebra $g(A)$ with connected Coxeter graph satisfy $0\leq \deg P(A)\leq 1$.

This conjecture is false for a finite type Lie algebra since in this case $\deg P(A)=$ the complex dimension of $F(A)$.

The infinite dimensional Kac-Moody Lie algebras contains affine type and indefinite type. For affine type Conjecture 3 is true since in this case $\deg P(A)=0$.

We need the following lemma for further discussion.

\noindent{\bf Lemma 3:} For rational functions $\displaystyle {P_i(t)},1\leq i\leq k$, $\deg P_i(t)=d_i$, suppose the sum of $P_i(t),1\leq t\leq k$ be $P(t)$, then

1. $\deg P(t)\leq \max\limits_{1\leq i\leq k} d_i$.

2. If $I=\{i|d_i=\max\limits_{1\leq t\leq k} d_t\}$ has only one element, then  $\deg P(t)= \max\limits_{1\leq i\leq k} d_i$.

By using Lemma 3, we have

\noindent {\bf Proposition 3: } The Poincar\'{e} series of an infinite dimensional Kac-Moody Lie algebra $g(A)$ with connected Coxeter graph satisfy $\deg P(A)\geq 0$.

\noindent {\bf Proof: } let $A$ be a $n\times n$ Cartan matrix of affine or indefinite type. For $I\subsetneqq S$, $P_I(A)$ denotes the Poincar\'{e} series of $A_I$. We prove this proposition by induction on
$n$. For $n=1$, the proposition is obviously true. Assume for each Cartan matrix $A'$ of affine or indefinite type with size $<n$, the proposition is true. By Equation (1.4), we have $$-\frac{(-1)^n}{P(A)}=\sum\limits_{I\subsetneqq S } \frac{(-1)^{|I|}}{P_I(A)} $$
\indent If $I=\emptyset$, $P_{\emptyset}(A)=1$; If $A_I$($I\not=\emptyset$) is of finite type, then $\deg P_I(A)\geq 1$. If $A_I$ is of affine or indefinite type, then by induction assumption $\deg P_I(A)\geq 0$. Therefor for any $I\subsetneqq S$, we have $\deg P_I(A)\geq 0$, so $\displaystyle{ \deg \frac{(-1)^{|I|}}{P_I(A)}\leq 0}$. By Lemma 3, $\displaystyle{ \deg \frac{(-1)^n}{P(A)}\leq 0}$, this shows $\deg {P(A)}\geq 0$.\qed

Let $A$ be a Cartan matrix and $I\subset S$, denote by $D_I$ the complex dimension of flag manifold $F(A_I)$. Put $I_0=S$. For $D=\infty$, set $t^{D}=0$.

By Equation (1.4)£¬
$$\frac{(t^{D_{I_0}}+(-1)^{{{|I_0|}}+1})}{P(A)}=\sum\limits_{I_1\subsetneqq S } \frac{(-1)^{|I_1|}}{P_{I_1}(A)} $$
We get

\noindent{\bf Lemma 4:} $\displaystyle{\frac{1}{P(A)}=\sum\limits_{S\supsetneqq {I_1}} \frac{(-1)^{|{I_1}|}}{(t^{D_{I_0}}+(-1)^{{|I_0|}+1})}\frac{1}{P_{I_1}(A)}}. $

By induction we get

\noindent{\bf Proposition 3:} 
{ $$
    \frac{1}{P(A)}=\sum\limits_{r=0}^{n-1} \sum\limits_{S\supsetneqq {I_1}\supsetneqq {I_2}\supsetneqq \cdots \supsetneqq {I_r} \supsetneqq \emptyset }
    (\prod\limits_{j=1}^r \frac{(-1)^{|{I_{j}}|}}{(t^{D_{I_{j-1}}}+(-1)^{{|I_{j-1}|}+1})})
\frac{1}{(t^{D_{I_r}}+(-1)^{|I_{r}|+1})}\hspace{1.2cm} (1.5)$$}

\indent Motivated by Equation (1.5), we give the following definition.

\noindent{\bf Definition 4:} Let $A$ be an $n\times n$ Cartan matrix of affine or indefinite type, $S=\{1,2,,\cdots,n\}$. A sequence $T: S=I_0\supsetneqq I_1\supsetneqq I_2\supsetneqq \cdots \supsetneqq I_r\supsetneqq \emptyset$ of subsets in $S$
is called a chain of length $r$.
If $\dim g(A_{I_i}),0\leq i\leq r$ is infinite, the chain $T$ is called an infinite chain. If $\dim g(A_{I_i}),0\leq i\leq r-1$ is infinite and $|I_r|=1$, the chain $T$ is called a quasi-infinite chain. The set of all the chains(infinite chains and quasi-infinite chains respectively) is denoted by $C(A)$($C_\infty(A)$ and $C^q_\infty(A)$ respectively).


Note $C_\infty(A)$ and $C_\infty^q(A)$ are disjoint.

Let $D$ be the dimension function which send $I\subset S$ to $\dim F(A_I)$. By Equation (1.5), we have

\noindent{\bf Proposition 4:} $1/P(A)=\sum\limits_{T\in C(A)} H_T$, where for $T: S=I_0\supsetneqq I_1\supsetneqq I_2\supsetneqq \cdots \supsetneqq I_r\supsetneqq \emptyset$, $$H_T=(\prod\limits_{j=1}^r \frac{(-1)^{|{I_{j}}|}}{(t^{D_{I_{j-1}}}+(-1)^{{|I_{j-1}|}+1})})
\frac{1}{(t^{D_{I_r}}+(-1)^{|I_{r}|+1})}.$$

This shows that the inverse power series $\displaystyle{\frac{1}{P(A)}}$ of Poincar\'{e} series $P(A)$ equals to the summation on the set of chains in $S$. The contribution of a chain is determined by the dimension function $D$.

\noindent{\bf Theorem 4:} The Conjecture 3 is false for a Cartan matrix $A$ if and only if for $A$ both $\sum\limits_{T\in C_{\infty}(A)} H_T$ and $\sum\limits^q_{T\in C_{\infty}(A)} H_T$ equal to $0$.

\noindent{\bf Proof:} The contribution of each chain $T: S=I_0\supsetneqq I_1\supsetneqq I_2\supsetneqq \cdots \supsetneqq I_k\supsetneqq \emptyset$ to $1/P(A)$ is $H_T$,
with degree $\leq 0$. Note $\deg H_T=0$ if and only if $T$ is a infinite chain and $\deg H_T=-1$ if and only if $T$ is a quasi-infinite chain.
By using the fact $H_T=(-1)^{-n-r-1}$ for infinite chain $T$ and Lemma 3 we know: if $\sum\limits_{T\in C_{\infty}(A)} H_T\not =0$, then $\deg P(A)=0$. If $\sum\limits_{T\in C_{\infty}(A)} H_T =0$, then the degree $0$ items in summation cancel each other. A quasi-infinite chain $T$ gives an item $$H_T=\frac{(-1)^{|{I_1}|}}{(-1)^{{|I_0|}+1}}
\frac{(-1)^{|{I_2}|}}{(-1)^{|I_1|+1}}\cdots
\frac{(-1)^{I_r}}{(-1)^{|I_{r-1}|+1}} \frac{{1}}{t+1}=\frac{(-1)^{1-n-r}}{(1+t)}.$$
if $$\sum\limits_{T\in C^q_{\infty}(A)} H_T =\sum\limits_{T\in C^q_{\infty}(A)}  \frac{(-1)^{1-n-r}}{(1+t)}\not=0. $$
By Lemma 3, we have $\displaystyle{\deg \frac{1}{P(A)}=-1}$, hence $\deg P(A)=1$. \qed

By the way we get the following criteria for $\deg P(A)=0$ or $1$. Let $K_0(A)=\sum\limits_{T\in C_{\infty}(A)} H_T$ and $K_1(A)=\sum\limits_{T\in C^q_{\infty}(A)} H_T$, we have

\noindent{\bf Proposition 5:} Let $A$ be a Cartan matrix of affine or indefinite type with connected Coxeter graph. If $K_0(A)\not=0$, the degree of the Poincar\'{e} series $P(A)$ is $0$, and the ratio of the coefficients of the highest power of $t$ in the numerator and denominator polynomials of $P(A)$ is ${1}:K_0(A)$; If $K_0(A) =0$, but $K_1(A)\not=0$, the degree of the Poincar\'{e} series $P(A)$ is $1$, and the ratio of the coefficients of the highest power of $t$ in the numerator and denominator polynomials of $P(A)$ is ${1}:K_1(A)$.

\noindent {\bf Remark: }Besides $K_0(A)$ and $K_1(A)$, we can also construct invariants $K_2(A),K_3(A),\cdots,$ $K_k(A),\cdots$ of $A$ as in obstruction theory, such that
$\deg P(A)=k$ if and only if $K_0(A),K_1(A),$ $\cdots,K_{k-1}(A)=0$,  but $K_{k}(A)\not=0$ and
the ratio of the coefficients of the highest power of $t$ in the numerator and denominator polynomials of $P(A)$ is ${1}:K_k(A)$. But things
become more and more complicate as $k$ increases, so we give up the explicit computation of $K_k(A)$ for $k>1$.

\noindent{\bf Example 6: } For rank $3$ Cartan matrix $A$, suppose the multiplicities of edges between three vertices be $p_{1},p_2,p_3$, then the only infinite chains are $S=\{1,2,3\}$ and $S\supset S-\{i\},p_i\geq 4$, hence $K_0(A)=0$ if and only if there is exact one pair of vertices connected by edge with multiplicity $\geq 4$.
In this case, there are $3$ length $1$ quasi-infinite chains $S\supsetneqq \{i\}, 1\leq i\leq 3$ and $2$ length $2$ quasi-infinite chains, so $K_1(A)=-1$. Therefore the Conjecture 3 is true for $n=3$. The ratio of the coefficients of the highest power of $t$ in the numerator and denominator polynomials of $P(A)$ can be computed from $K_0(A)$ and $K_1(A)$.

For $n=4$, there are counter examples to the Conjecture 3.

\noindent{\bf Example 7: } For Cartan matrices

$$A=\left(\begin{array}{rrrr}
                2 & -1 & 0 & -1 \\
                -4 & 2 & -1 & 0 \\
                0 & -1 & 2 & -1 \\
                -1 & 0 & -2 & 2
              \end{array}\right) \text{and }B=\left(
                  \begin{array}{rrrr}
                    2 & -1 & -1 & 0 \\
                    -1 & 2 & -1 & -1 \\
                    -1 & -4 & 2 & -1 \\
                    0 & -2 & -2 & 2 \\
                  \end{array}
                \right)
$$
Their Poincar\'{e} series are
$$P(A)=-\frac{(t^5+t^4+t^3+t^2+t+1)(t^3+t^2+t+1)(t+1)}{t^5+t^4+t^3+t^2+t-1}.$$ and
$$P(B)=-\frac{(t^5+t^4+t^3+t^2+t+1)(t^3+t^2+t+1)(t+1)}{t^7+t^6+2t^5+2t^4+2t^3+2t^2+t-1}.$$
with degree $4$ and $2$.

We don't have examples to show that the degree of Poincar\'{e} series can be arbitrary large and don't know if there exists any upper bound.





\end{document}